\documentclass{amsart}
\usepackage[english]{babel}
\usepackage{amsmath, layout}
\usepackage{amssymb}
\usepackage{latexsym}
\usepackage{amscd}
\newtheorem{theorem}{Theorem}[section]         
\newtheorem{lemma}[theorem]{Lemma}             
\newtheorem{corollary}[theorem]{Corollary}     
\newtheorem{proposition}[theorem]{Proposition} 
\newtheorem{definition}[theorem]{Definition}   
\newtheorem{example}[theorem]{Example}         
\numberwithin{equation}{section}

\newenvironment{namelist}[1]{%
\begin{list}{}
    {
      
      \settowidth{\labelwidth}{#1}
      \setlength{\leftmargin}{1.1\labelwidth}
    }
  }{%
\end{list}}

\newcommand{\Z}{\ensuremath{\mathbb{Z}}}

\newcommand{\C}{\ensuremath{\mathfrak {C}}}
\newcommand{\A}{\ensuremath{\mathcal {A}}}
\newcommand{\E}{\ensuremath{\mathcal {E}}}

\newcommand{\rk}{\operatorname{rk}}
\newcommand{\OS}{\operatorname{OS}}
\newcommand{\M}{\ensuremath{\mathcal M}}
\newcommand{\I}{\ensuremath{\Im}}
\newcommand{\cl}{\ensuremath{c\ell}}
\begin{document}
\title[Quadratic Orlik-Solomon  algebras]
{Quadratic Orlik-Solomon  algebras\\ of graphic matroids}
\thanks{2000
\emph{Mathematics Subject Classification}: \emph{Primary}:\ 05B35, 14F40; \emph{Secondary}:\
52C35.\\
 \emph{Keywords and phrases}: arrangement of
 hyperplanes, 
cohomology algebra,  matroid, 
Orlik-Solomon algebra.
}
\author{Raul Cordovil and David Forge}
\address{{}\newline
Departamento de Matem\'atica,\newline 
Instituto Superior T\' ecnico \newline
 Av.~Rovisco Pais
 - 1049-001 Lisboa  - Portugal}
\email{cordovil@math.ist.utl.pt}
\thanks{The  first author's research  was 
supported in part by FCT (Portugal) through program POCTI,
 the projects SAPIENS/36563/99 and
  by DONET (European Union Contract
ERBFMRX-CT98-0202). 
The second  author's research  was  supported  by FCT through the projects SAPIENS/36563/99
 and
  by DONET}
\address{{}\newline
Laboratoire de recherche en informatique\newline 
Batiment 490
Universite Paris Sud\newline 
91405 Orsay Cedex -
France
}
\email{forge@lri.fr}

\begin{abstract}
{In this note we introduce a 
sufficient condition for the Orlik-Solomon algebra associated to a matroid $\M$ to be $l$-adic and 
we prove that this condition is necessary when
$\M$ is binary (in particular graphic). Moreover, this result
 cannot be extended to the   class of all matroids. 
}
\end{abstract}

\maketitle 

\section{Introduction}
Throughout this note  $\M$  denotes
 a simple matroid of rank $r$ on the ground set
$[n]$.   We refer to 
\cite{W2}  
 as a standard source  for  matroids. 
An  {\em (essential) arrangement of hyperplanes} in $\mathbb{K}^d$
is a finite collection $\A_{\mathbb{K}}= \{H_1,\ldots ,H_n\}$ of codimension 1
subspaces of $\mathbb{K}^d$ such that $\bigcap_{H_i\in \A_{\mathbb{K}}} H_i= 0.$ 
There is a  
   matroid $\M(\A_{\mathbb{K}}),$ on the ground set 
  $[n]:=\{1,\dotsc, n\}$,
canonically  determined by 
$\A_{\mathbb{K}}$: i.e.,
a subset  $X\subseteq [n]$ is independent if and only if
the codimension of $\bigcap_{i\in X}H_i$ is equal to the cardinality of $X$.
The manifolds $\mathfrak{M}(\A_{\mathbb{C}})=
\mathbb{C}^d\setminus\bigcup_{H_i\in
 \A} H_i$  are important in the Aomoto-Gelfand
theory of $\A_{\mathbb{C}}$-hypergeometric functions, see \cite{OT2} for a recent introduction from the
point of   view of arrangement theory. We refer  to the survey \cite{FR} for a
recent  
discussion on the role of matroid theory in the study of hyperplane 
arrangements.
Fix a set $E:=\{e_1,\dotsc, e_n\}$ and consider the Grassmann
algebra
$\mathcal{E}:=\bigwedge\big(\bigoplus_{i=1}^n\mathbb{K}e_i\big).$ 
We set
$\mathcal{E}_\ell:=\bigwedge^\ell\big(\bigoplus_{e\in E}\mathbb {K}e\big),\,\, \forall \ell\in \mathbb{N}.$ 
(By convention
$\mathcal{E}_0=\mathbb{K}$\, even in the case where $E=\emptyset.$)
For every linearly ordered subset $X=({i_1},\dotsc,  {i_m})\subset [n], ~
i_1<\dotsm < {i_m},$\, let $e_X$ be the monomial
$e_X:=e_{i_1}\wedge e_{i_2}\wedge \dotsm\wedge e_{i_m}.$ By definition 
  set  $e_\emptyset=1\in \mathbb{K}.$ 
For $\mathfrak{X}$  a  subset of\, $2^{[n]},\,$ let\,
$\I(\mathfrak{X}):=\big\langle\partial(e_X): X
\in
\mathfrak{X}\big\rangle$\, be the two-side graded ideal of the Grassmann
algebra $\E\,$ generated by  the ``differentials'',
$$\partial(e_{_X})=\partial(e_{i_1}\wedge \dotsm \wedge
e_{i_m})=\sum(-1)^{j}e_{i_1}\wedge\dotsm \wedge
{e}_{i_{j-1}}\wedge
{e}_{i_{j+1}}\dotsm
\wedge e_{i_m}.$$ 
The graded Orlik-Solomon
$\mathbb{K}$-algebra
$\OS(\M)$   
is the quotient   
$\mathcal{E} 
/\Im(\C),$ where $\C$ denotes the set of circuits of $\M.$
The de Rham cohomology algebra $H^{\bullet}\big(\mathfrak{M}
(\mathcal{A}_{\mathbb{C}}); \mathbb{C}\big)$  is
shown to be isomorphic to the  Orlik-Solomon $\mathbb{C}$-algebra of the matroid 
$\M(\A_{\mathbb{C}}),$ see \cite{OS1}. 
Since then the combinatorial structure of the manifold 
$\mathfrak{M}(\mathcal{A}_{\mathbb{C}})$  have been intensively studied
(i.e., the properties depending of the matroid  $\M(\A_{\mathbb{C}})$).  
Quadratic OS (Orlik-Solomon) algebras (i.e., such that $\I(\C)=\I(\C_3)$ where $\C_\ell,\, \ell\in [n],$ denotes the
set of circuits with at most $\ell$ elements of $\M$)
appear in the study of complex  arrangements, in the rational homotopy theory of the manifold
$\mathfrak{M}$ and the Koszul property of $\OS(\M),$ 
see \cite{F4,PY}.
It is an old open problem to find a condition on the matroid $\M\,$ equivalent to
$\OS(\M)$ being quadratic, or in general {\em $\ell$-adic} \big(i.e., such that
$\I(\C)=\I(\C_{\ell+1})\big)$.
 Some  necessary conditions are known but none of these  conditions implies the
quadratic property. 
For more details, we refer the reader to the recent survey on  OS algebras \cite{Yuz}, which also provides a large
bibliography.

We say that a circuit $C$  of a  matroid has a {\em chord}\,  $i_\alpha\in C,$ 
if there are two circuits $C_1$ and $C_2$
such that $C_1\cap C_2=i_\alpha$ and $C=C_1\Delta C_2.$
We say that a  matroid is $\ell$-\emph{chordal}, $\ell \geq 4,$ if  
every circuit with at least $\ell$ elements has a chord. $\M$ is said \emph{chordal} if
it is 4-chordal.
 We prove that if  
all the circuits of $\M$ with  at least
$\ell+2$ elements have a chord, then the algebra $\OS(\M)$ is $\ell$-adic. 
The $2$-adic case with $\rk(\M)=3$ is known and correspond to the
``parallel 3-arrangements", see \cite{Falk5} for details. A matroid $\M\,$ is said to be \emph{binary} if it is realizable over $\Z_2,$ or
equivalently  if the symmetric difference of two different 
 circuits of $\M\,$ contains a circuit. Graphic (and cographic) matroids
are  important examples of binary matroids.
We prove also that in the case of binary  matroids   
the reverse is also true.
 The proofs use only simple algebraic arguments and basic matroid theory.
\section{$\ell$-adic  OS algebras}
Let $G=(V; S)$ be a connected simple graph and suppose the vertices  [resp. edges] 
of G labelled with the integers
$1,\dotsc,d,$ [resp. $1,\dotsc,n,$]
i.e., $V=[d]$\, and
$S=[n].$ Consider the hyperplane arrangement in $\mathbb{K}^d$ 
 given by $$\A_{\mathbb{K}}:=\{\text{Ker}(x_i-x_j):~ \{i,j\} \in S\}.$$
We say that $\A_{\mathbb{K}}$ is the \emph{graphic arrangement} determined by $G.$
It is clear that the  matroid $\M(\A_{\mathbb{K}})$ on the ground set $[n]$ 
is the cycle matroid of the graph $G.$   We recall some of our notations.
Let $\mathfrak{X}$ be a  subset of\, $2^{[n]}.\,$ Let\,
$\I(\mathfrak{X})=\big\langle\partial(e_X): X
\in\mathfrak{X}\big\rangle$ be the two-side ideal 
of the Grassmann
$\mathbb{K}$-algebra $\E$ generated by the elements $\{\partial(e_X): X \in
\mathfrak{X}\}.$ We recall that the differential of degree $-1,\,$ 
$\partial : \E \to \E,$ can be defined by the equalities
\begin{equation}\label{a}
\partial(e_i)=1 \hspace{3mm} \text{for every} \hspace{3mm} e_i,\, i=1,\dotsc, n
\end{equation}
 and the \emph{Leibniz rule}
\begin{equation}\label{b}
\partial(a\wedge b)=\partial(a)\wedge b+ (-1)^{\text{deg}(a)}a\wedge \partial(b)
\end{equation}
for every pair of homogeneous elements $a,b \in \E.$ 
From  Equations~ (\ref{a}) and (\ref{b}) we conclude that
\begin{equation}\label{c}
\partial^2(e)=0\hspace{3mm} \text{
for every element}\hspace{3mm} e\in \E.
\end{equation}
We make use of the following technical (but fundamental) lemma.
\begin{lemma}\label{bord}
Consider two  subsets $X$ and $X'$ of\, $[n]$ of at least two elements
and an element $i_{\alpha}\in X.$
Let $\mathfrak{X}$ be a  subset of\, $2^{[n]}.\,$
Then the following properties hold:
\begin{namelist}{xxxxxxxxx}
\item[$(\ref{bord}.1)$] $\partial(e_X)\in \I(\mathfrak{X})\iff e_X\in \I(\mathfrak{X}).$
\vspace{1mm}
\item[$(\ref{bord}.2)$]  $\big(e_{X\setminus i_\ell} \in \I(\mathfrak{X}),\,\, 
\forall i_\ell\in X\setminus i_{\alpha} \big)
\,
\Longrightarrow\,  e_{X\setminus i_{\alpha}}\in \I(\mathfrak{X}).$
\vspace{1mm}
\item[$(\ref{bord}.3)$] If\,
$X\cap X'=i_{\alpha},$ then\, $e_{{X\Delta X'}}\in \I(\{X,X'\}).$ 
\end{namelist}
 \end{lemma}
\begin{proof}
$(\ref{bord}.1)$\, Note that
$$\partial(e_X)\in \I(\mathfrak{X})\,\,\Longrightarrow\,\, e_{i_{\alpha}}\wedge
\partial(e_X)=\pm e_X\in \I(\mathfrak{X}).$$ 
Conversely suppose that
$$e_X=\sum_i y_i\wedge 
\partial(e_{Y_i}),\,\,\text{where}\, \, y_i\in \E,\,
\text{and}\,\,{Y_i}\in \mathfrak{X}.$$
From  Equations (\ref{a}), (\ref{b}) and (\ref{c}) we conclude that
$$\partial(e_X)=\sum_i\partial (y_i)\wedge\partial(e_{Y_i})\in \I(\mathfrak{X}).$$

\vspace{1mm}
\noindent{$(\ref{bord}.2)$}\,
Choose an element $i_\ell\in X\setminus i_{\alpha}.$ 
By hypothesis we know that  $e_{X\setminus i_\ell} \in \I(\mathfrak{X})\,$
and so\, $e_X= \pm e_{i_{\ell}} \wedge e_{X\setminus i_\ell} \in \I(\mathfrak{X}).$
From $(\ref{bord}.1)$ we know that $\partial (e_X)\in \I(\mathfrak{X})$ 
 and so $$e_{X\setminus i_{\alpha}}=\pm \partial(e_{X})+\sum_{i_{\ell}\not =i_{\alpha}}
\pm e_{X\setminus i_\ell}
\in \I(\mathfrak{X}).$$
$(\ref{bord}.3)$\, We know that $|X\cup X'|\geq 3.$ It is clear that 
 $e_{X\cup X'\setminus i_j}\in \I(\{X, X'\}),$  for all $i_j\not =i_{\alpha}.$ 
From $(\ref{bord}.2)$  we conclude that
 $e_{X\cup X'\setminus i_{\alpha}}\in \I(\{X, X'\}).$ By hypothesis 
we know that $X\cap X'= i_{\alpha}$\, so\,
$X\Delta X'=X\cup X'\setminus i_{\alpha}.$
\end{proof}
\begin{definition}\label{Db}
{\em
A set\, $\mathfrak{X}\subseteq 2^{[n]}$  is said to be $\Delta$-\emph{closed} if, 
for all subsets $X,X'\subseteq {[n]},$ we have:
\begin{namelist}{xxxxxxx}
\item[~(\ref{Db}.1)]
$\big(X\subseteq X'\,\,\, \text{and}\,\,\, X\in \mathfrak{X}\big)\,\Longrightarrow\,
X'\in\mathfrak{X},$
\item[~(\ref{Db}.2)] $(X\cap X'=i_{\alpha}, \,\,|X|\geq 2,\,|X'|\geq 2\,\, 
\text{and}\,\,\, X,X'\in \mathfrak{X})\,\Longrightarrow\, X\Delta X'
\in\mathfrak{X}.$
\end{namelist}
Given a subset  $\mathfrak{Z}\subseteq 2^{[n]}$ we will note $\cl_{\Delta}(\mathfrak{Z})$
the  smallest $\Delta$-closed subset of $2^{[n]}$  containing $\mathfrak{Z}.$
}
\end{definition}

\begin{proposition}\label{bb}
$\M$ is $\ell$-chordal then $\,\C\subseteq\cl_\Delta(\C_{\ell-1}).$
\end{proposition}
\begin{proof}
As every circuit with $\ell$ or more elements has a chord the result is clear by induction on the number of
elements of the circuits.
\end{proof}
As a consequence of  Lemma~\ref{bord} and  Definition~\ref{Db} we get:
\begin{proposition}\label{s}
For every subset\,  $\C'\subseteq \C(\M),$  we have\, $\I(\C')=\I(\cl_{\Delta}(\C')).$\qed
\end{proposition}
The following Corollary extends  Proposition 3.2 in \cite{Falk5} (case $\ell=4$ 
and rk$(\M)=3).$ It is a consequence of Propositions \ref{Db} and \ref{bb}.
\begin{corollary}\label{apart} If the matroid $\M$ is $\ell$-chordal then the associated\, 
$\OS$ algebra $\OS(\M)$ is $(\ell-2)$-adic. \qed
\end{corollary}
\begin{proposition}\label{chord} Let \M\, be a 
binary matroid.  For every circuit $C$ of $\M$ with $k$ elements,  $\,k\geq 4,$   
the following statements are equivalent:
\begin{namelist}{xxxxxxxxx}
\item[$~(\ref{chord}.1)$] $C$\, has a chord, 
\item[$~(\ref{chord}.2)$] $\partial(e_C)\in
 \I(\C_{k-1}),$ 
\item[$~(\ref{chord}.3)$] $e_C\in 
 \I(\C_{k-1}).$ 
\end{namelist}
\end{proposition}
\begin{proof} $~(\ref{chord}.2) \iff (\ref{chord}.3)$ is a consequence of $(\ref{bord}.1).$

\noindent{$(\ref{chord}.1) \Longrightarrow (\ref{chord}.3).\,$} Suppose that  
$C$ has a chord $i_\alpha$ and let $C_1, C_2$ be two  circuits
such that $C_1\cap C_2= i_\alpha$ and $C_1\Delta C_2= C.$  From the equivalence 
$(\ref{bord}.1)\iff (\ref{bord}.2)$ we know that $e_{C_1}, e_{C_2} \in
\I(\C_{k-1}).$ It follows from
$(\ref{bord}.3)$ that $e_{C\setminus i_\alpha}\in\I(\C_{k-1}).$\newline
$(\ref{chord}.3) \Longrightarrow (\ref{chord}.1).$  We will prove it by
contradiction. Let
$C$  be a circuit without a chord and suppose  that $e_C\in
\I(\C_{k-1}).$  From the definitions we know that
$$e_C=\sum_j  a_j\wedge
\partial(e_{C_j}),\,\, a_j\in \E,\,\, C_j \in \C_{k-1}.$$
 For this to be possible, there must exist a circuit
$C_{j_\ell}$ such that $|C\cap C_{j_\ell}|=|C_{j_\ell}|-1.$ Set $x=C_{j_\ell}\setminus
C.$ As \M\,
is  binary we know that  $C'=C\Delta C_{j_\ell}=C\setminus  C_{j_\ell}\cup x $ 
is also a circuit. We conclude that $C'\cap  C_{j_\ell}= x$
and $C=C'\Delta C_{j_\ell}$ and so $x$ is a chord of $C,$ a contradiction.
\end{proof}
\begin{corollary}\label{binary}
Let \M\, be a
 binary matroid and $\ell, \, \ell \geq 4,$ a natural number. Then the following
statements are equivalent:
\begin{namelist}{xxxxxxx}
\item[$(\ref{binary}.1)$] The matroid $\M$ is $\ell$-chordal.
\item[$(\ref{binary}.2)$] The  algebra\,\, $\OS(\M)$ is $(\ell-2)$-adic.
\end{namelist}
In particular \M\, is chordal iff the algebra\,\, $\OS(\M)$ is quadratic.
\end{corollary}
\begin{proof} $(\ref{binary}.1) \Longrightarrow (\ref{binary}.2)$\, 
is a special case of Corollary~\ref{apart}.\newline
$(\ref{binary}.2) \Longrightarrow (\ref{binary}.1)$ From the definitions 
we know that for all $k,\, k\geq \ell,$ 
$\I(\C_{k-1})=\I(\C_{\ell-1})=\I(\C).$ We conclude  from the Proposition
\ref{chord} that every circuit with $k,\, k\geq \ell,$ elements has a chord, 
i.e., $\M$ is $\ell$-chordal.
\end{proof}
Corollary \ref{apart} gives a  condition which makes  a matroid quadratic and we have
shown in Corollary \ref{binary} that this condition is also necessary in 
the case of binary matroids. It is known that this condition  is not
necessary in general, see \cite{FR}.  
Close related with our results we propose a weaker sufficient condition for quadraticity.
\begin{definition}\label{Db2}
{\em
A set\, $\mathfrak{X}\subseteq 2^{[n]}$  is said to be $\Delta'$-\emph{closed} if, 
for all subsets $X,X'\subseteq {[n]}$ and $i_\alpha\in [n]$ we have:
\begin{namelist}{xxxxx}
\item[~$(S_1)$]
$\big(X\subseteq X'\,\,\, \text{and}\,\,\, X\in \mathfrak{X}\big)\,\Longrightarrow\,
X'\in\mathfrak{X},$
\item[~$(S_3)$] $(X\setminus i_\ell\in \mathfrak{X},  \,\,\forall i_\ell \in X\setminus i_\alpha)\,\Longrightarrow\,
X\setminus i_\alpha\in \mathfrak{X}.$
\end{namelist}
Given a subset  $\mathfrak{Z}\subseteq 2^{[n]}$ we will note $\cl_{\Delta'}(\mathfrak{Z})$
the  smallest $\Delta'$-closed subset of $2^{[n]}$  containing $\mathfrak{Z}.$
}
\end{definition}
\begin{proposition}\label{mm}
For every subset $\C'\subseteq \C(\M),$ we have $\cl_{\Delta}(\C') \subseteq \cl_{\Delta'}(\C')$ and
$\I(\C')=\I(\cl_{\Delta'}(\C')).$ In particular if $\C\subseteq \cl_{\Delta'}(\C_3)$ then the algebra $\OS(\M)$ is
quadratic.
\end{proposition}
\begin{proof}
Remark that Conditions (\ref{Db}.1) and (\ref{Db}.3) imply Condition (\ref{Db}.2). So the result is a consequence of
Lemma ~\ref{bord} and Definition~\ref{Db2}.
\end{proof}
The following example is a case where $\C\not\subseteq \cl_{\Delta}(\C_3)$ but $\C\subseteq \cl_{\Delta'}(\C_3).$
So $\OS(\M)$ is quadratic from Proposition~\ref{mm}.
We know no example where $\C\not \subseteq \cl_{\Delta'}(\C_3)$ and $\OS(\M)$ is quadratic.
\begin{example}
{\em
Consider the rank 3 (simple) matroid $\M$ of Figure 1  
 on 7 elements whose circuits of 3 elements are 123, 145, 167, 246 and 357.
  It is easy to see that $\C\subseteq \cl_{\Delta'}(\C_3)$ (i.e., $\OS(\M)$ is quadratic) but $\C\not \subseteq
\cl_{\Delta}(\C_3)$ (i.e., 
 the matroid  is  not
chordal). For example we can check that the circuit 2356 is not in $\cl_{\Delta}(\C_3)$ but is in
$\cl_{\Delta'}(\C_3).$ }
\end{example}
\vspace{4mm}
\begin{center}
\begin{picture}(-53,90)(100,30)

\put(10,10) {\line(1,0){120}}
\put(10,10) {\circle*{6}}\put(-1,10){3}
\put(70,10) {\circle*{6}}\put(68,-3){5}
\put(130,10) {\circle*{6}}\put(135,10){7}

\put(10,10) {\line(3,4){60}}
\put(40,50) {\circle*{6}}\put(29,50){2}
\put(70,90) {\circle*{6}}\put(67,96){1}

\put(40,50){\line(1,0){60}}
\put(70,50) {\circle*{6}}\put(58,38){4}

\put(70,10){\line(0,1){80}}
\put(70,90) {\line(3,-4){60}}
\put(100,50) {\circle*{6}}\put(105,50){6}
\put(70,10) {}\put(50,-23){\textbf{Figure 1}}
\end{picture}
\end{center}
\vspace*{3cm}

\centerline{ACKNOWLEDGEMENTS}

\vspace{0.4cm}
The first author is grateful to Marcelo Viana for  an invitation to visit the IMPA and for providing financial
support.
He also thanks  the ``Organizing Committee" of the  ``Latin-American Workshop on Cliques in Graphs" for an invitation to speak at the
conference.

\end{document}